\documentclass[a4paper,12pt]{article}
\usepackage[utf8]{inputenc}
\usepackage[T2A]{fontenc}  
\usepackage[english]{babel}
\usepackage{hyperref}
\usepackage{graphicx}

\title{On periodic approximate solutions of dynamical systems with a quadratic right-hand side}
\author{Mikhail Malykh and Leonid Sevastianov}
\usepackage[centertags]{amsmath}
\usepackage{amsfonts}
\usepackage{amssymb}
\usepackage{amsthm}
\usepackage{xcolor}
\newcommand{\alert}[1]{{\color{red}{#1}}}
\renewcommand{\alert}[1]{#1}

\theoremstyle{remark}

\newtheorem*{acknowledgments}{Acknowledgments}
\theoremstyle{definition}

\newtheorem{example}{Example}

\def\dt{{\Delta t}}

\begin{document}\large\maketitle

\begin{abstract}
Difference schemes are considered for dynamical systems $ \dot x = f (x) $ with a quadratic right-hand side, which have $t$-symmetry and are reversible. Reversibility is interpreted in the sense that the Cremona transformation is performed at each step in the calculations using a difference scheme. The inheritance of periodicity and the Painlevé property by the approximate solution is investigated. In the computer algebra system Sage, such values  are found for the step $ \dt $, for which the approximate solution is a sequence of points with the period $ n \in \mathbb N $. Examples are given and hypotheses  about the structure of the sets of initial data generating sequences with the period $ n $ are formulated.  

\end{abstract}

\section{Introduction}

At the end of the last century, a study of the inheritance  of algebraic properties of exact solutions of ordinary differential equations, including dynamical systems, by difference schemes was begun. The     inheritance of algebraic integrals is the most studied question. Back in the 1990s, a family of symplectic Runge-Kutta schemes was discovered that preserve all linear and quadratic integrals of motion, e.g., all classical integrals in the problem of the rotation of a top
\cite{Cooper,Suris-1988,Suris-1990,Lubich}. 
At the same time, the first difference scheme  that preserves all algebraic integrals of the many-body problem was found
by Greenspan  \cite{Greenspan-1992,Greenspan-1995-1,Greenspan-1995-2,Greenspan-2004}
and independently by  Simo and Gonz{\'a}lez
\cite{Simo-1993,Graham}.
Other conservative schemes for this problem, including high-order ones,  can be constructed by introducing additional variables, with respect to which all integrals of the many-body problem can be written as quadratic ones \cite{malykh-casc-2020}, 
thus combining the  invariant energy quadratization method,  proposed
recently  by Yang et al.  
\cite {Yang-2016, Yang-2017,Zhang-2020}
and the scalar auxiliary variable (SAV) approach, proposed by Jie Shen et al. 
\cite{Shen}.
Unfortunately, this way leads to implicit schemes, the use of which for numerical calculations turns out to be very resource-consuming.
\cite{Zhang-2020}.

Much less studied is the question of the inheritance of reversibility by approximate solutions. In classical mechanics, it is believed that there is a one-to-one correspondence between the initial and final positions of bodies. In fact, an exact solution to a general nonlinear dynamical system has this property only locally. Moreover, as Painlevé noted
\cite{Painleve,Umemura,malykh-pomi-2015},
the possibility of integrating a dynamical system in Abelian functions and solution reversibility are closely related to each other. 

Any difference scheme defines a correspondence between the values of the solution at two times, separated by the step $ \dt $. Therefore, it is natural to call it reversible if it defines a birational correspondence between these values (the definition is given below in Section \ref{n:def}). When considering birational correspondences, it is convenient to consider the phase space as projective. At the same time, we noticed that any dynamical system with a quadratic right-hand side is approximated by an invertible difference scheme (Section \ref{n:hat}).  

The theory of birational transformations of the plane was laid by L. Cremona
\cite{Cremona},
therefore, birational transformations on projective spaces are called Cremona transformations.
It should be noted that even  two-dimensional Cremona transformations   are a very complicated object. In the theory of dynamical chaos,  simple quadratic transformations, e.g., the integer Hénon transformation
\cite{Henon-1973,Tabor}
are used to explain the origin of chaos. To relate this discrete model to Hamiltonian systems, Tabor
\cite[\S 4.2]{Tabor}
noticed that models of this kind arise when discretizing differential equations with respect to $ t $. However, the difference scheme described by him preserves the symplectic structure and therefore inherits quadratic integrals, but is not reversible. These properties cannot be combined at the same time
\cite{malykh-casc-2019}, so we intend to sacrifice simplecticity in favor of reversibility. This property is extremely important not only  for mechanics, but also for creating efficient numerical methods for studying dynamical systems, since it is devoid of the main drawback of conservative schemes, their implicitness.

Among dynamical systems with a quadratic right-hand side, there are models describing oscillations of a pendulum and rotation of a top and integrable in elliptic functions. The main qualitative property of this model is the periodicity of the solution. In this paper, we want to use these simple examples to understand how exactly the periodicity can be inherited by an approximate solution found using a reversible scheme. 

\section{Definitions}
\label{n:def}

Let us consider  dynamical system 
\begin{equation}
\label{eq:ode}
\frac{dx_i}{dt} =f_i(x_1, \dots, x_m), \quad i =1,2, \dots, m.
\end{equation}
with  polynomial right-hand side
\[
f_i \in \mathbb{Q}[x_1, \dots, x_m].
\]
For brevity, we will use vector notation, meaning by $ x $ the tuple  $ (x_1, \dots, x_m) $. Within the framework of the finite difference method
\cite{Lubich}
the system of differential equations is replaced with the system of algebraic equations
\begin{equation}
\label{eq:hat}
g_i(x, \hat x, \dt)=0, \quad i=1,\dots, n.
\end{equation}
In this case, $ x $ is interpreted as the value of the solution at the time $ t $, and $ \hat x $  as the solution at the time $ t + \dt $.

From the point of view of mathematical modeling of mechanical phenomena,   the inheritance of two properties by the difference scheme, namely, $ t $-symmetry and reversibility, is of  greatest impor\-tance. We say that the difference scheme \eqref{eq:hat} has $ t $-symmetry if it is invariant under the transformation
\[
\dt \to -\dt, \quad x \to \hat x, \quad \hat x \to x.
\]
By reversibility, we should understand the possibility to uniquely determine the final data $ \hat x $ from the initial data $ x $  and vice versa using the  system \eqref{eq:hat}  for any fixed value of the step $ \dt $. Since Eqs. \eqref{eq:hat} are algebraic, this means that $ \hat x $ must be a rational function of $ x $, and $ x $  must be a rational function of $ \hat x $. We will consider $ x, \hat x $ as two points of the projective space $ \mathbb {P} _m $ and say that the difference scheme \eqref{eq:hat} is invertible if for any fixed value of $ \dt $ this scheme defines a  Cremona transformation .
The combination of $ t $-symmetry and reversibility means that the difference scheme defines a one-parameter family of Cremona transformations $\mathcal C$, such that  $\hat x=\mathcal  C(\dt) x$ and $\mathcal  C(\dt)^{-1}=\mathcal  C(-\dt)$.

In mechanics, the property of reversibility has to be introduced in a more complicated way, while the Cauchy theorem makes it possible to substantiate reversibility locally, in a vicinity of a nonsingular point. However, globally a dynamical system may not possess this property. Indeed, let us consider the following Painlevé
\cite{Painleve}
 initial problem
\begin{equation}
\label{eq:Cauchy}
\frac{dx}{dt}=f(x), \quad x|_{t=t_0}=x_0
\end{equation}
on the segment $[t_0, t_0 + \dt]$ of the real axis $t$. For some values of  $t_0$, the procedure for analytic continuation of the solution obtained in the Cauchy theorem along a segment does not encounter singular points other than poles, and in this case the final value of $x(t_0 + \dt)$ is uniquely determined by the initial value of $x_0$. However, if the path encounters a branch point, then the final value depends on the way it is passed. Therefore, $ x (t_0 + \dt) $ is, generally speaking, a multivalued function of the initial value $ x_0 $. Thus, if a dynamical system has the global reversibility property, then it also has the Painlevé property
\cite[\S 3.5]{Goriely}:
the singular  points of the solution are not branch points.

Classical completely integrable models, including pendulums and tops, are integrable in elliptic functions and, as can be seen from the solution, have the Painlevé property
\cite{Golubev}.
As a rule, the general solution to completely integrable nonlinear models possessing the Painlevé property defines a birational transformation not on the entire phase space, but on integral manifolds distinguished by algebraic integrals.

\begin{example}
\label{ex:wp:1}
The dynamical system
\begin{equation}
\label{eq:wp}
\dot x= y, \quad \dot y= 6x^2-a
\end{equation}
has an algebraic integral 
\begin{equation}
\label{eq:wp:c1}
\frac{y^2}{2}- 4x^3 +ax =C_1,
\end{equation}
having the meaning of total mechanical energy. Therefore, this system can be integrated in the Weierstrass elliptic functions
\[
x = \wp(t +C_2 ,2a,C_1), \quad y= \wp'(t+C_2,2a,C_1),
\]
where $ C_1, C_2 $ are the integration constants, and for this reason  below it is referred to as the $ \wp $-oscillator. By virtue of the addition theorem for elliptic functions, the general solution can be expressed rationally in terms of  
\[
\wp(t,2a,C_1), \, \wp'(t,2a,C_1) \quad  \mbox{and} \quad \wp(C_2,2a,C_1), \, \wp'(C_2,2a,C_1).
\] 
Therefore, the general solution is expressed rationally in terms of the initial data $ x_0 = \wp (C_2,2a, C_1) $, $ y_0 = \wp'(C_2,2a, C_1) $, and the coefficients of this expression depend on $ t $ and $ C_1 $  transcendentally. It follows from the $ t $-symmetry that these expressions define a birational transformation on the elliptic curve \eqref{eq:wp:c1}. Thus, on the integral curve \eqref{eq:wp:c1}, the Cauchy problem defines a birational correspondence between initial and final data. However, this corres\-pon\-dence does not continue until the Cremona transformation of the entire $ xy $ plane, although at first it seemed quite surprising
\cite[ch. 7]{Klein-1}. 
%Это дает при фиксированном $C_1$ однопараметрическое семейство решений, обладающих одним и тем же периодом. Величина $C_1$ определяет величину этого периода и зависимость общего решения от $C_1$ оказывается существенно трансцендентной. 
\end{example}

\section{Equations with  quadratic right-hand side}
\label{n:hat}

At present, there are well developed  Painlevé tests, which make it possible to find out in practice whether a given dynamical system has the Painlevé property. Generally speaking, these tests are a set of necessary conditions for the absence of moving branch points, which  are checked algorithmically
\cite[\S 3.9]{Goriely}.

In the one-dimensional case ($ n = 1 $), only the Riccati equation has the Painlevé property
\begin{equation}
\label{eq:Riccati}
\frac{dx}{dt}=a+bx+cx^2
\end{equation}
for any, including zero values of constants $ a, b, c $. Moreover, the initial problem defines a M\"obius transformation on the projective line. It is not difficult to construct a difference scheme that inherits this property:
\begin{equation}
\label{eq:Riccati:hat}
\hat x -x=\left(a+b\frac{x+\hat x}{2}+cx\hat x\right)\dt
\end{equation}
Since any birational transformation on a projective line is a M\"obius one, it is easy to prove the converse: in the one-dimensional case, an invertible difference scheme can be constructed only for the Riccati equation
\cite{malykh-pomi-2018}.

However, for $ n> 1 $, the continuous and discrete case lose their similarity. For any dynamical system with a quadratic right-hand side, a $ t $ -symmetric reversible difference scheme can be constructed:
\begin{equation}
\label{eq:quadr}
\hat x_i -x_i=F_i(x,\hat x) \dt, \quad i=1,\dots, n,
\end{equation}
where $ F_i $ is obtained from $ f_i $ by replacing monomials: $ x_j $ with $ (\hat x_j + x_j) / 2 $, $ x_jx_k $ with $ (\hat x_j + x_j) ( \hat x_k + x_k) / 4 $, and $ x_j ^ 2 $ with $ x_j \hat x_j $. However, only a few dynamical systems with a quadratic right-hand side possess the Painlevé property. This issue was studied a long time ago: the dynamical system describing the rotation of a rigid body around a fixed point always has a quadratic right-hand side and has the Painlevé property only in 3 special cases found by S.V.~Kovalevskaya
\cite{Golubev}.

\begin{figure}
\begin{center}
\includegraphics[width=0.5\textwidth]{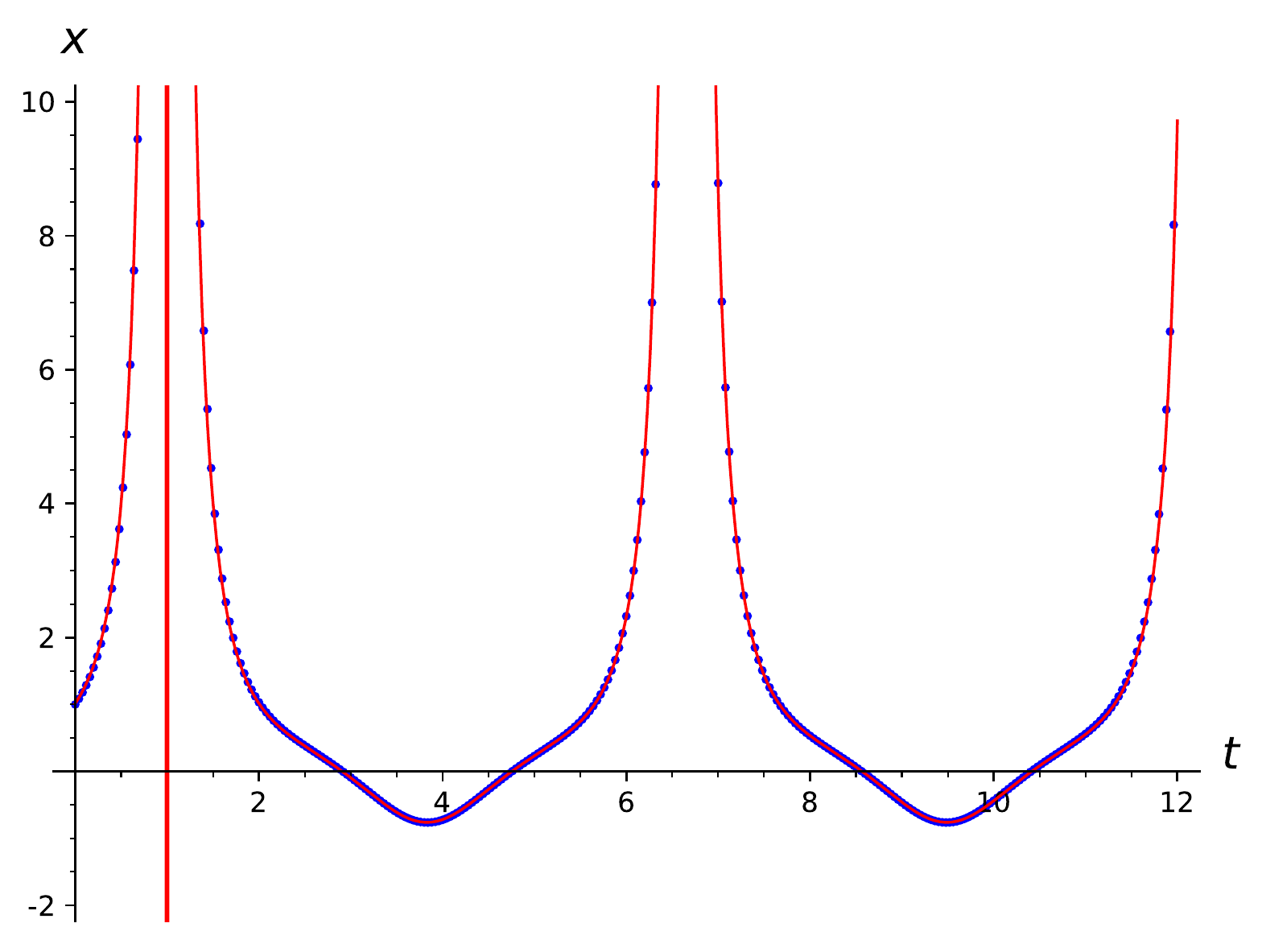}
\end{center}
\caption{Solution of system \eqref{eq:wp} at $a=1/2$ and initial conditions $(x,y)=(1,2)$. The continuous line plots the exact solution, the dots represent the approximate one. }
\label{fig:wp:1}
\end{figure}

We investigated how the Painlevé property is inherited by an approximate solution using two examples, the Riccati equation \eqref{eq:Riccati} and the $ \wp $-oscillator \eqref{eq:wp}. In both cases, it turned out that the calculations using a reversible  difference scheme can be continued through a pole without noticeable accumulation of errors. For the Riccati equation, this property was substantiated in
\cite{malykh-pomi-2018}. It seems to us that this property of an approximate solution successfully transfers the concept of the Painlevé property to finite difference equations, which by no means excludes other interpretations
\cite{Grammaticos,Clarkson,Ishizaki}.

\begin{example}
\label{ex:wp:2}
Figure \ref{fig:wp:1} shows the solution of  system \eqref{eq:wp} for $ a = 1/2 $ and initial conditions $ (x, y) = (1,2) $. The calculation of an approximate solution according to a reversible scheme does not encounter any difficulties over the entire considered interval $ 0 <t <12 $, containing two poles of the exact solution. 
\end{example}

\section{Periodicity}

The approximate solution is a sequence $ x_0, x_1, \dots $, each next element of which is obtained from the previous one by applying the Cremona transformation $\mathcal{C}$:
\[
x_{n+1}=\mathcal C x_{n}
\] 
This sequence will have period  $n$, if $x_{n}=x_0$, i.e., if $x_0$ is a fixed point of  $\mathcal C^n$. 

Let a positive integer $ n $ and an initial value $ x_0 \in \mathbb Q ^ m $ be given. Then the step $ \dt $ at which the sequence has a period $ n $ can be selected in the following way. Considering $ \dt $ as a symbolic variable, we calculate $ \mathcal C ^ n x_0 $. We get $ m $ rational functions from $ \mathbb Q (\dt) $. Equating them to $ x_0 $, we obtain $ m $ of algebraic equations, the common roots of which are the required step values. Generally speaking, several equations for one variable may not have common roots.

We have considered three examples: i) a linear oscillator that can be easily investigated analytically 
\cite{mal-2019}, 
ii) $\wp$-oscillator, and iii) Jacobi oscillator, i.e., dynamical system 
\begin{equation}\label{eq:pqr}
\dot p = qr, \, \dot q =-pr, \, \dot r =-k^2 pq,
\end{equation}
integrable in terms of elliptic Jacobi functions. We chose different initial data and considered  $n$ in the interval from $2$ to $10$. All calculations were performed in the Sage computer algebra system 
\cite{sage}
on an office PC. The degrees of polynomials, the common roots of which give the desired step values, increase exponentially with $ n $, which significantly limited our ability to increase $ n $. 

\begin{table}
\centering
\begin{tabular}{ |l|l| } 
 \hline
 $n$ & $\dt$  \\ 
 \hline
 2 & $\emptyset$  \\ 
 3 & $\emptyset$  \\ 
 4 & $1.074$  \\ 
 5 & $6.908$  \\ 
 6 & $\emptyset$  \\ 
 7 & $0.556, 5.870, 7.759$  \\ 
 8 & $0.535, 1.074, 6.843$  \\ 
 9 & $0.504, 9.187$  \\ 
10 & $0.471, 0.559, 6.777, 6.908$  \\ 
 \hline
\end{tabular}
\caption{The values of  step $ \dt $ providing  periodicity of the solution from  example \ref{ex:wp:3}}
\label{t:wp}
\end{table}

\begin{example}
\label{ex:wp:3}
For the $ \wp $-oscillator under the same initial conditions as were used in  example \ref{ex:wp:2}, there are no values of step $ \dt $  for which the solution has a period $ n = 2,3, 6 $. For $ n = 4 $ the step is independent of the starting point. Table \ref{t:wp} contains all the matched positive values for $ \dt $ found for $ n $ in the top ten.
\end{example}

\begin{table}
\centering
\begin{tabular}{ |l|l| } 
 \hline
 $n$ & $\dt$  \\ 
 \hline
 2 & $\emptyset$  \\ 
 3 & $3.609$  \\ 
 4 & $2.041$  \\ 
 5 & $1.47$, $6.86$  \\ 
 6 & $1.17$, $3.60$  \\ 
 7 & $0.97$, $2.57$, $10.85$  \\ 
 8 & $0.83$, $2.04,$ $5.18$  \\ 
 9 & $0.73$, $1.70$, $3.60$, $16.23$  \\ 
 \hline
\end{tabular}
\caption{The values of  step $ \dt $ providing  periodicity of the solution from  Example \ref{ex:pqr:1}}
\label{t:pqr:1}
\end{table}

\begin{example}
\label{ex:pqr:1}
For the  Jacobi oscillator \eqref{eq:pqr} with $k = \tfrac{1}{5}$ under the initial conditions
\[
p=0, \quad q=1, \quad r=0,
\]
there are positive values of step $ \dt $  for any periods $ n>2$, see Table \ref{t:pqr:1}. 

\begin{figure}
\begin{center}
\includegraphics[width=0.3\textwidth]{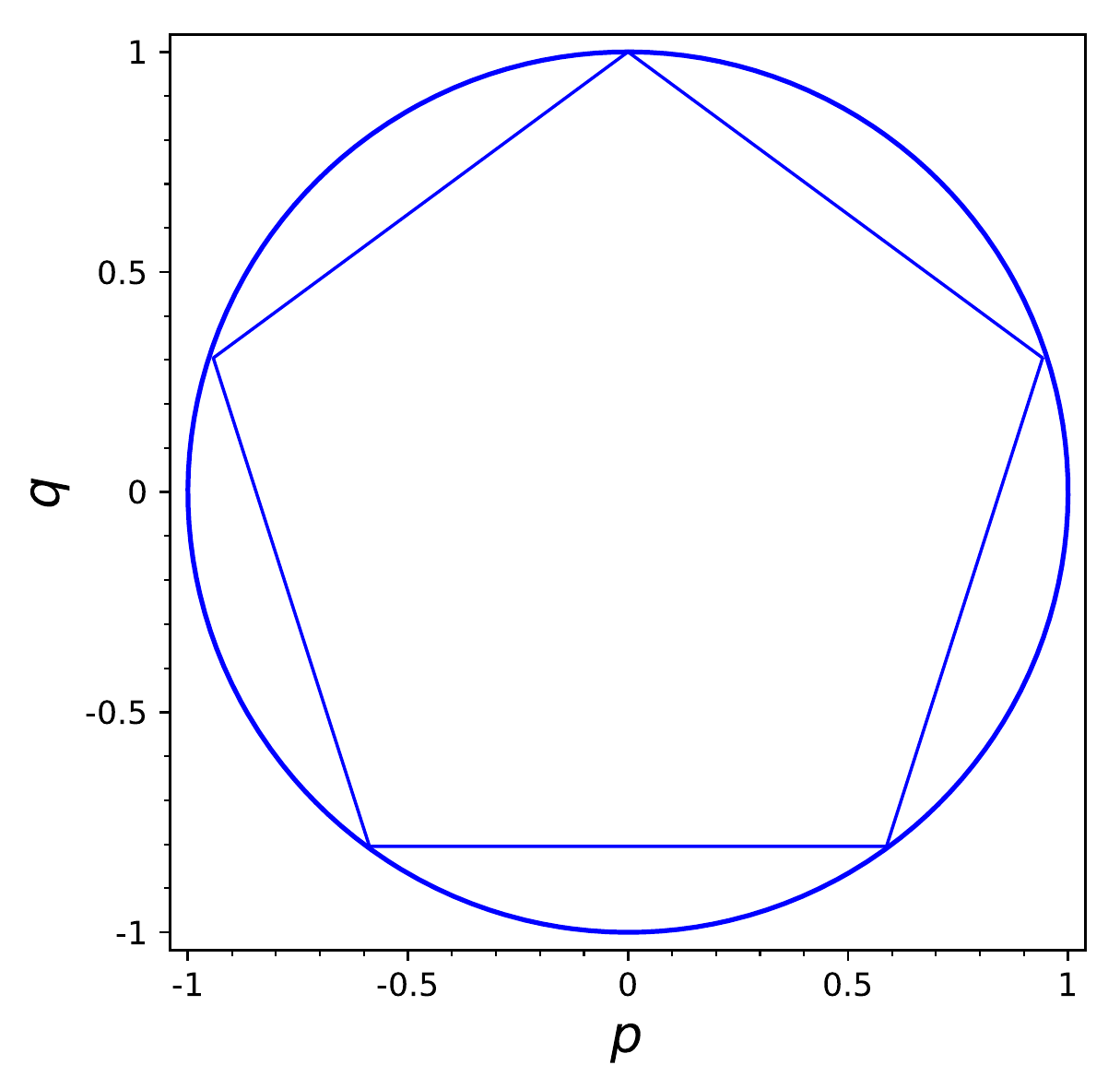}
\includegraphics[width=0.3\textwidth]{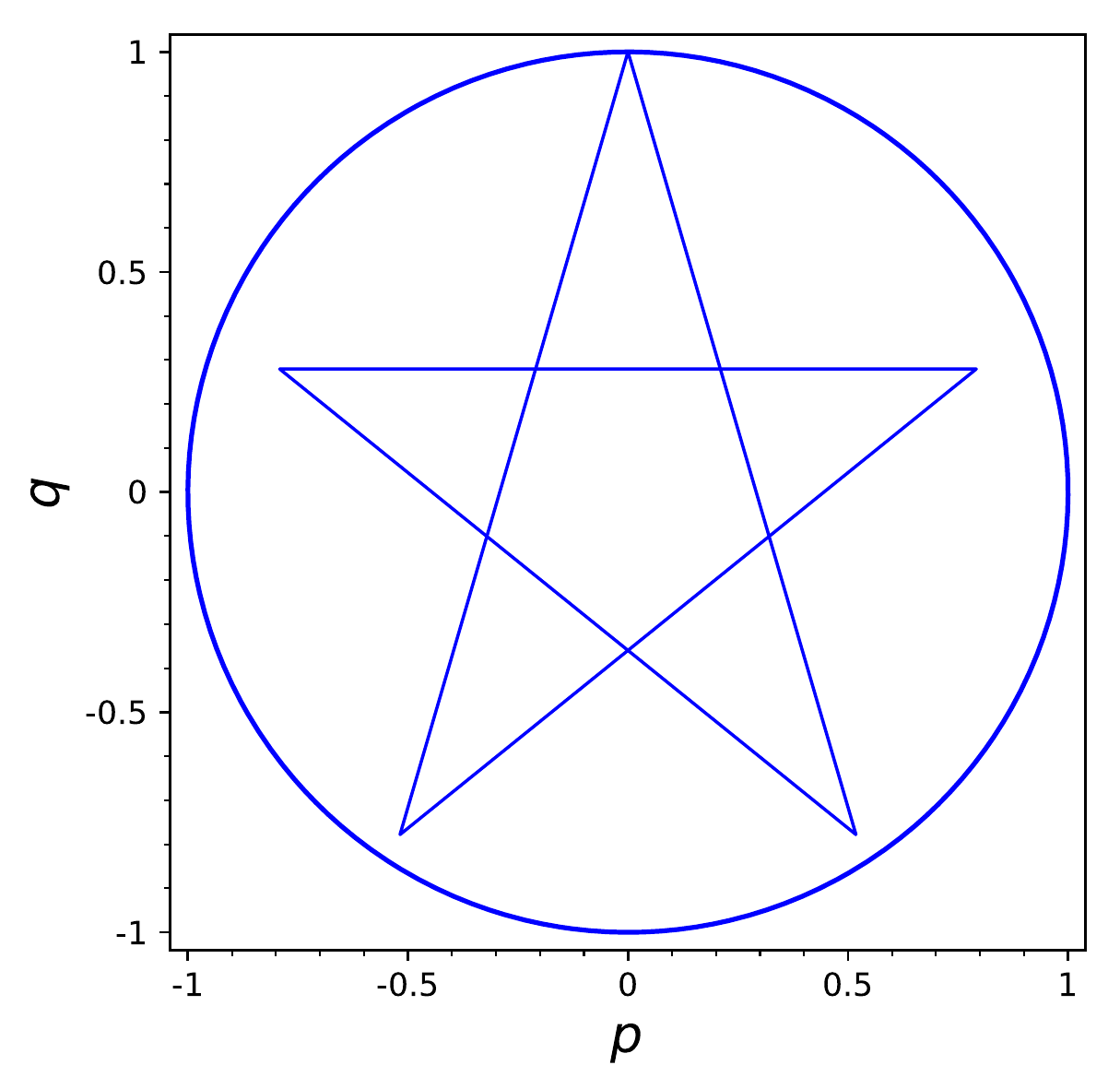}
\end{center}
\caption{Approximate solution from Example \ref{ex:pqr:1} at two steps ensuring the period $n=5$.}
\label{fig:pqr:5}
\end{figure}

\begin{figure}
\begin{center}
\includegraphics[width=0.3\textwidth]{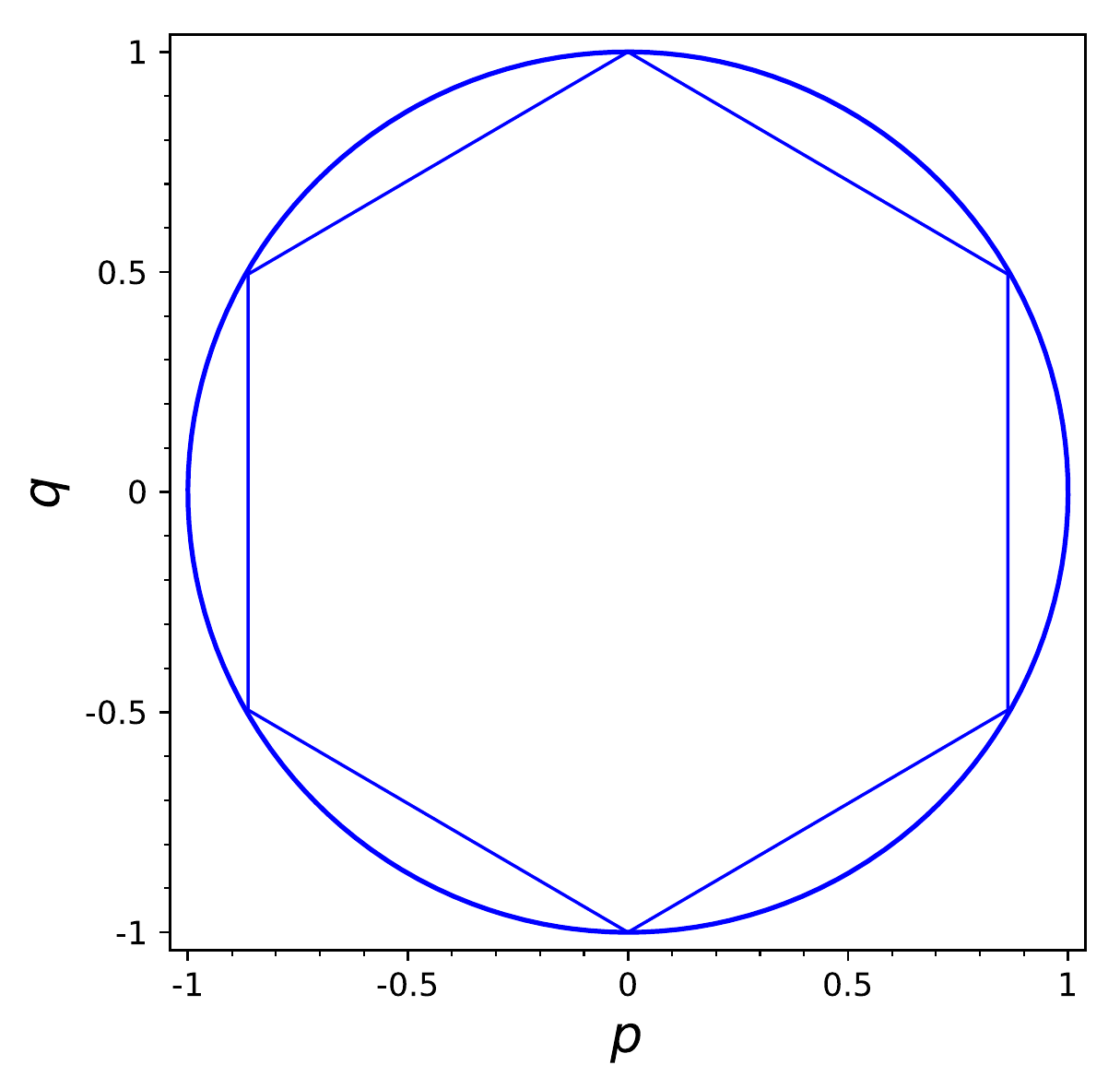}
\includegraphics[width=0.3\textwidth]{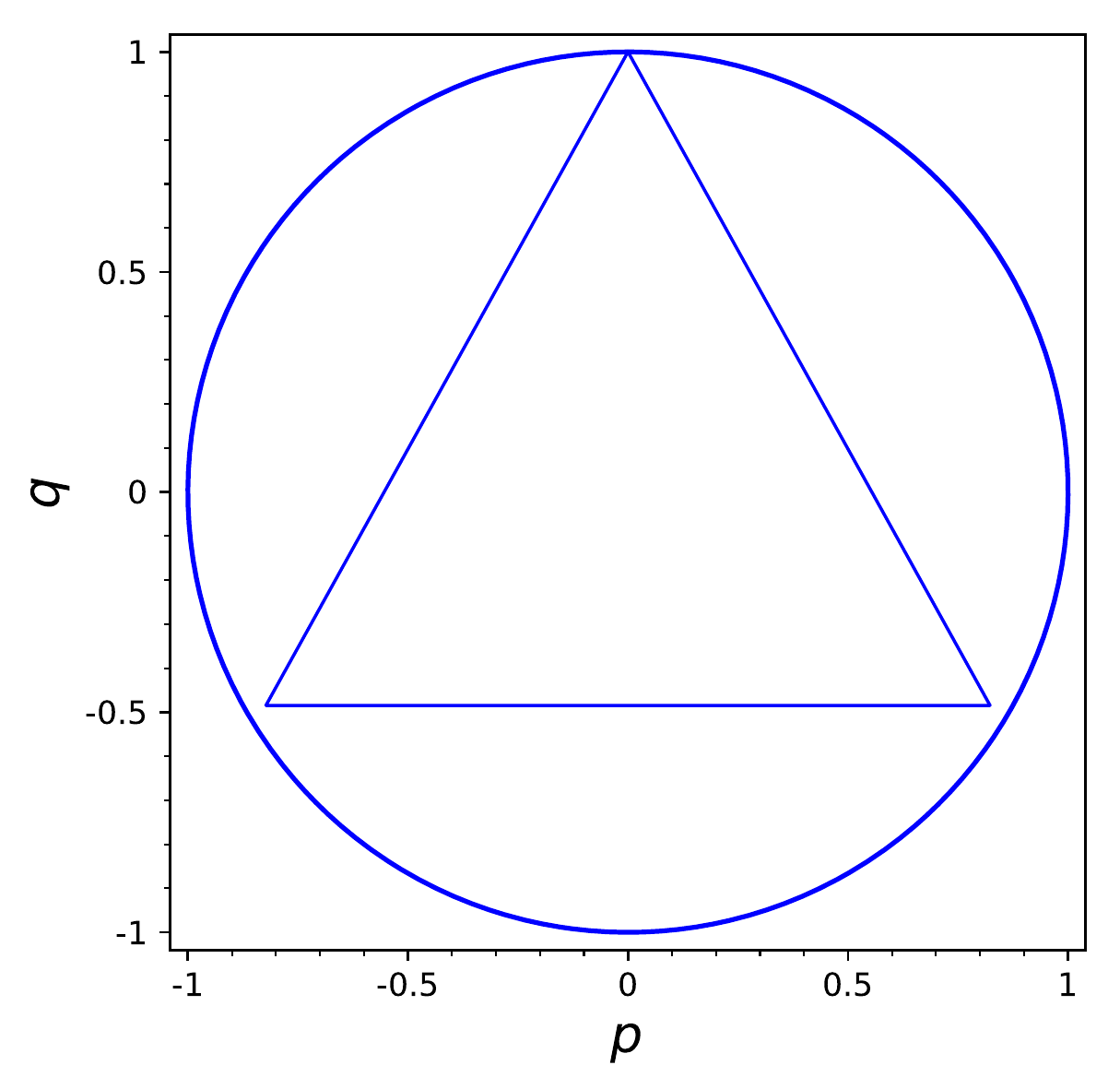}
\end{center}
\caption{Approximate solution from Example \ref{ex:pqr:1} at two steps ensuring the period $n=6$.}
\label{fig:pqr:6}
\end{figure}

\begin{figure}
\begin{center}
\includegraphics[width=0.3\textwidth]{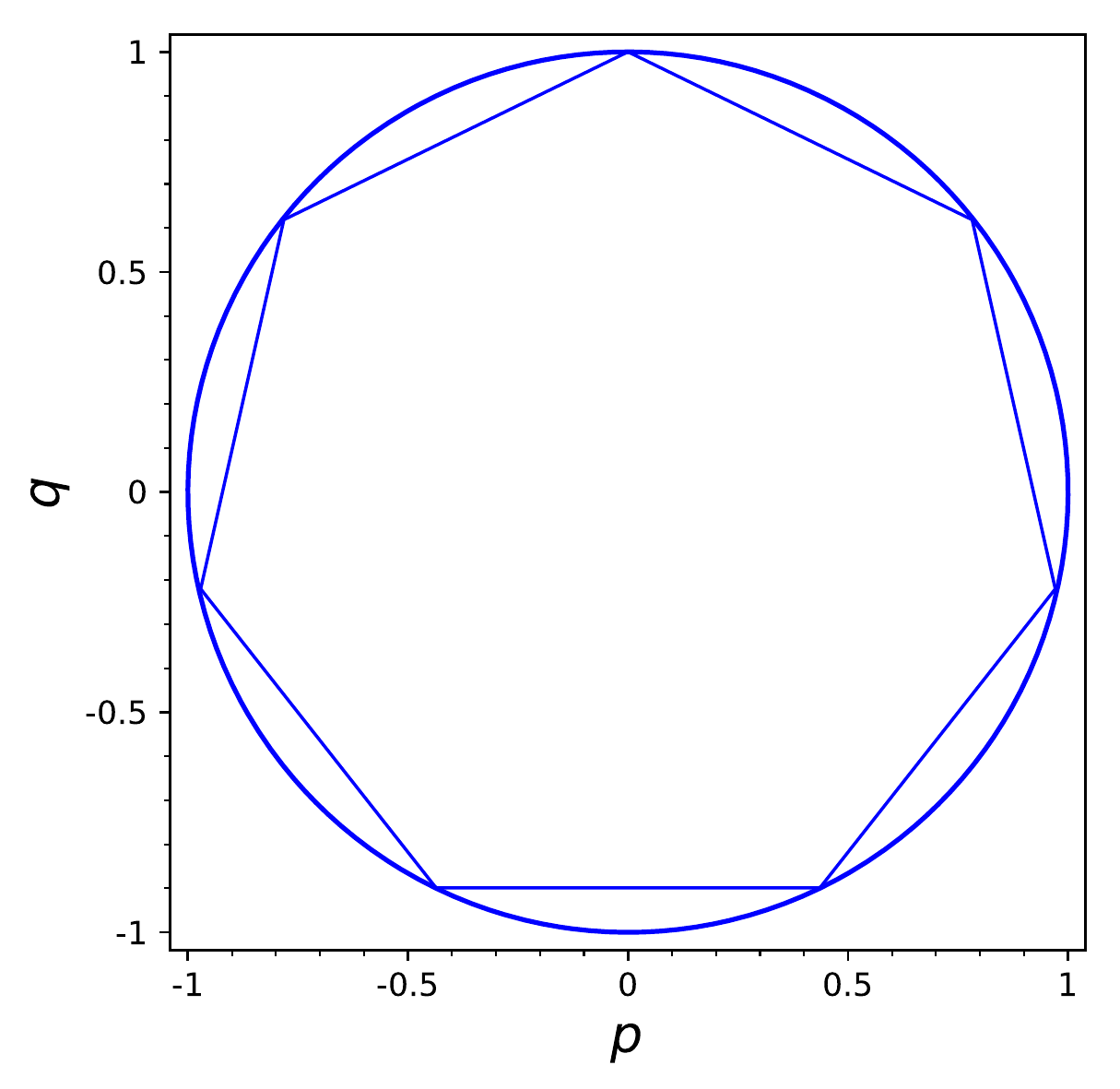}
\includegraphics[width=0.3\textwidth]{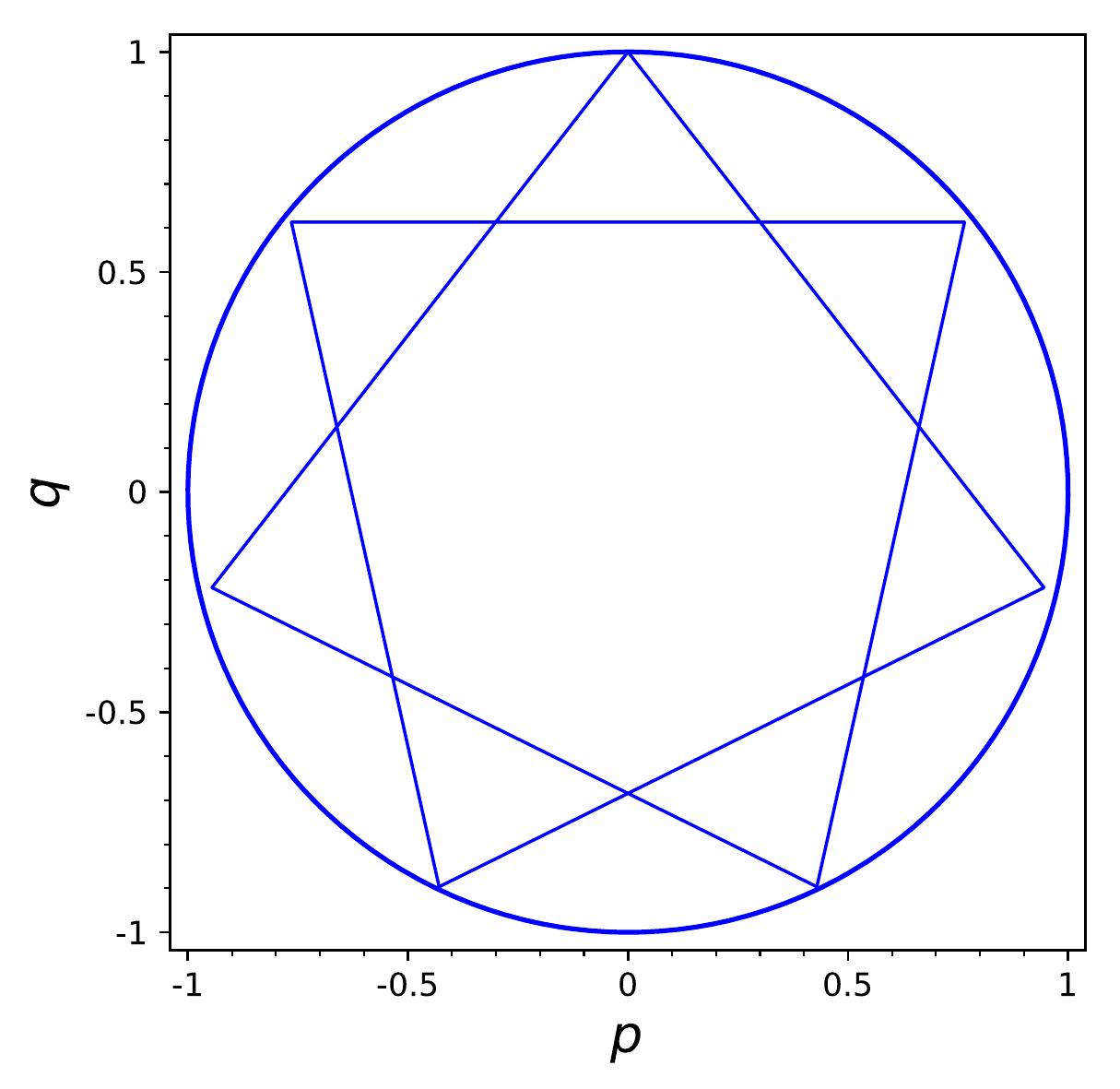}
\includegraphics[width=0.3\textwidth]{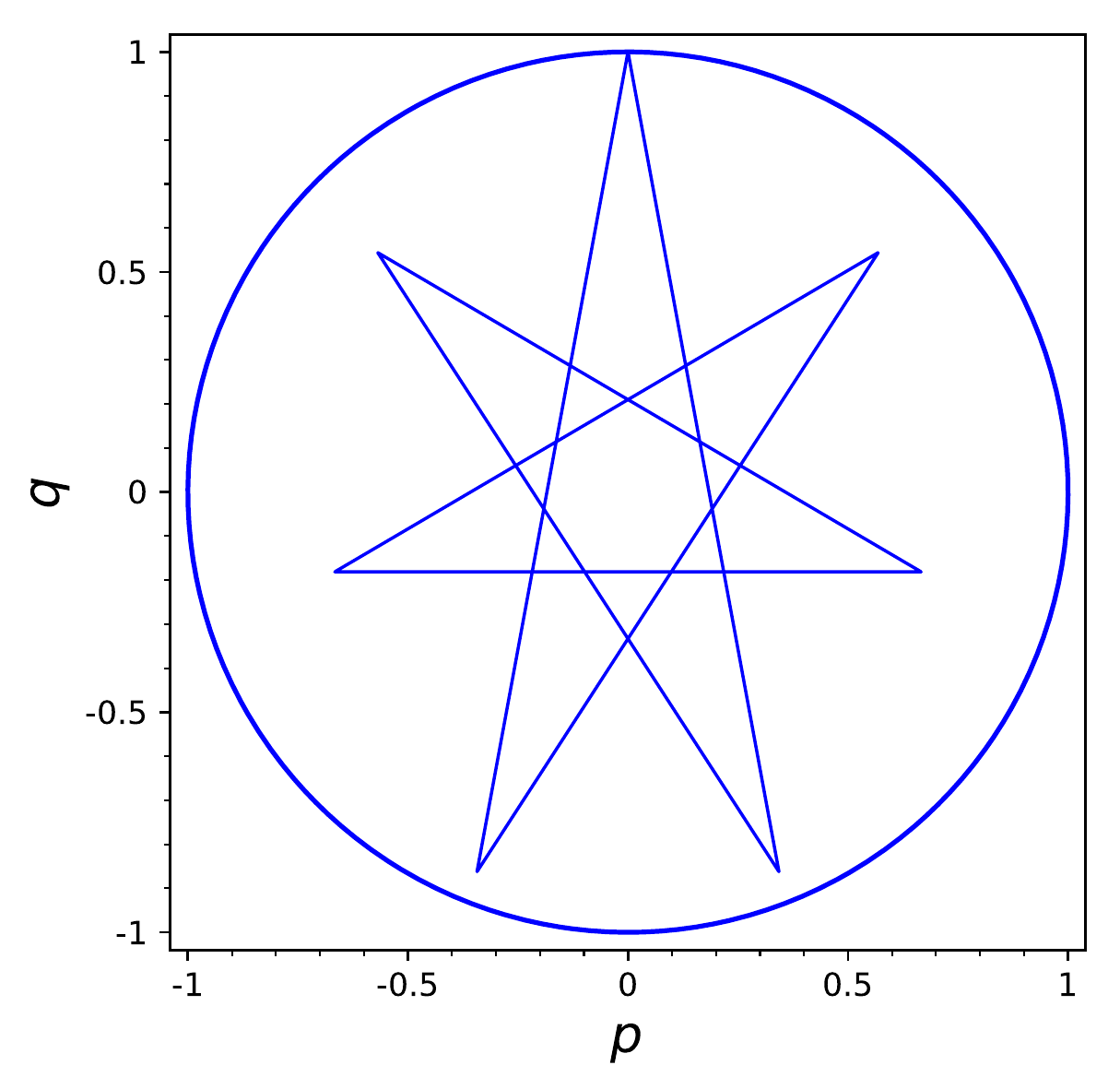}
\end{center}
\caption{Approximate solution from Example \ref{ex:pqr:1} at two steps ensuring the period $n=7$.}
\label{fig:pqr:7}
\end{figure}

Approximate solution of system \eqref{eq:pqr} for $k = 1/5$, satisfying initial conditions $(p,q,r)=(0,1,1)$, has the period $n=5$ at two values. 
In the plane $ pq $, at the first value of the step, an almost regular pentagon is obtained and at the second value we obtain a pentagram (Fig. \ref{fig:pqr:5}). In both cases, the integral of motion $ p ^ 2 + q ^ 2 = 1 $  is not exactly conserved. \alert{Approximate solution has the period $n=6$ at two values of the step. First of them is coincide with the step at $n=3$  and give us a triangle, the second give a hexagon. Approximate solution has the period $n=7$ at three values of the step. }

\begin{table}
\centering
\begin{tabular}{ |l|l| } 
 \hline
 $n$ & $n\dt$  \\ 
 \hline
 3 & $10.827$  \\ 
 4 & $8.164$  \\ 
 5 & $7.379$  \\ 
 6 & $7.022$  \\ 
 7 & $6.827$  \\ 
 8 & $6.706$  \\ 
 9 & $6.627$  \\ 
$\infty$ & $6.347$  \\ 
 \hline
\end{tabular}
\caption{Sequence of transitions for approximate solutions from Example \ref{ex:pqr:1}}
\label{t:pqr:2}
\end{table}
As $ n $ grows, the number of step values at which periodic approximations to the solution of the Cauchy problem are obtained, considered in the example \ref{ex:pqr:1}, grows. In this case, the smallest possible $ \dt $ for fixed $ n $ corresponds to an almost regular $ n $-gon in the $ pq $ plane. These solutions revert to their original value in times $ n \dt $, collected in Table \eqref{t:pqr:2}.  These times seem to form a monotonically decreasing sequence converging to the exact period.
\end{example}

It is convenient to present the results of the experiments carried out in the form of two hypotheses: 
\begin{enumerate}
\item for any sufficiently large $ n $ and any initial conditions, one can specify a finite number of positive values for the step $ \dt $, at which periodic sequences with the period $ n $ are obtained,
\item if we associate each $ n $ with a minimum period, we get a sequence converging to the period of the exact solution for $n \to \infty$.
\end{enumerate}
By virtue of the first hypothesis, any exact particular solution can be approximated by an approximate solution that inherits the periodic nature of the exact solution, and by virtue of the second hypothesis the approximation step $ \dt $  can be taken arbitrarily small and, therefore, approach the exact solution with any given accuracy.

\section{Equiperiodic sets}

In the previous Section, we followed one solution, but changed $ \dt $. Let us now look at the behavior of solutions in the phase space, but for a fixed $ \dt $.

The set in the phase space formed by all the initial data generating approximate solutions with the same period $ n $ is algebraic; we will call it an equiperiodic set of the $ n $-th order. It is easy to deduce from the first hypothesis that equiperiodic sets of sufficiently large order are not empty and have codimension 1.

To find it in the previous algorithm, it is necessary to consider $ x_0 $ as a tuple of $m$ symbolic variables. We managed to find these sets only for small $ n $.

\begin{table}
\centering
\begin{tabular}{ |l|l| } 
 \hline
 $n$ & Degre of $F_n$  \\ 
 \hline
 4 & 0  \\ 
 5 & 3  \\ 
 6 & 3  \\ 
 7 & 6  \\ 
 8 & 6  \\ 
 9 & 9  \\ 
10 & 12  \\ 
 \hline
\end{tabular}
\caption{The degrees of equiperiodic curves for Example \ref{ex:Fn}}
\label{t:Fn}
\end{table}

\begin{figure}
\centering
\includegraphics[width=0.45\textwidth]{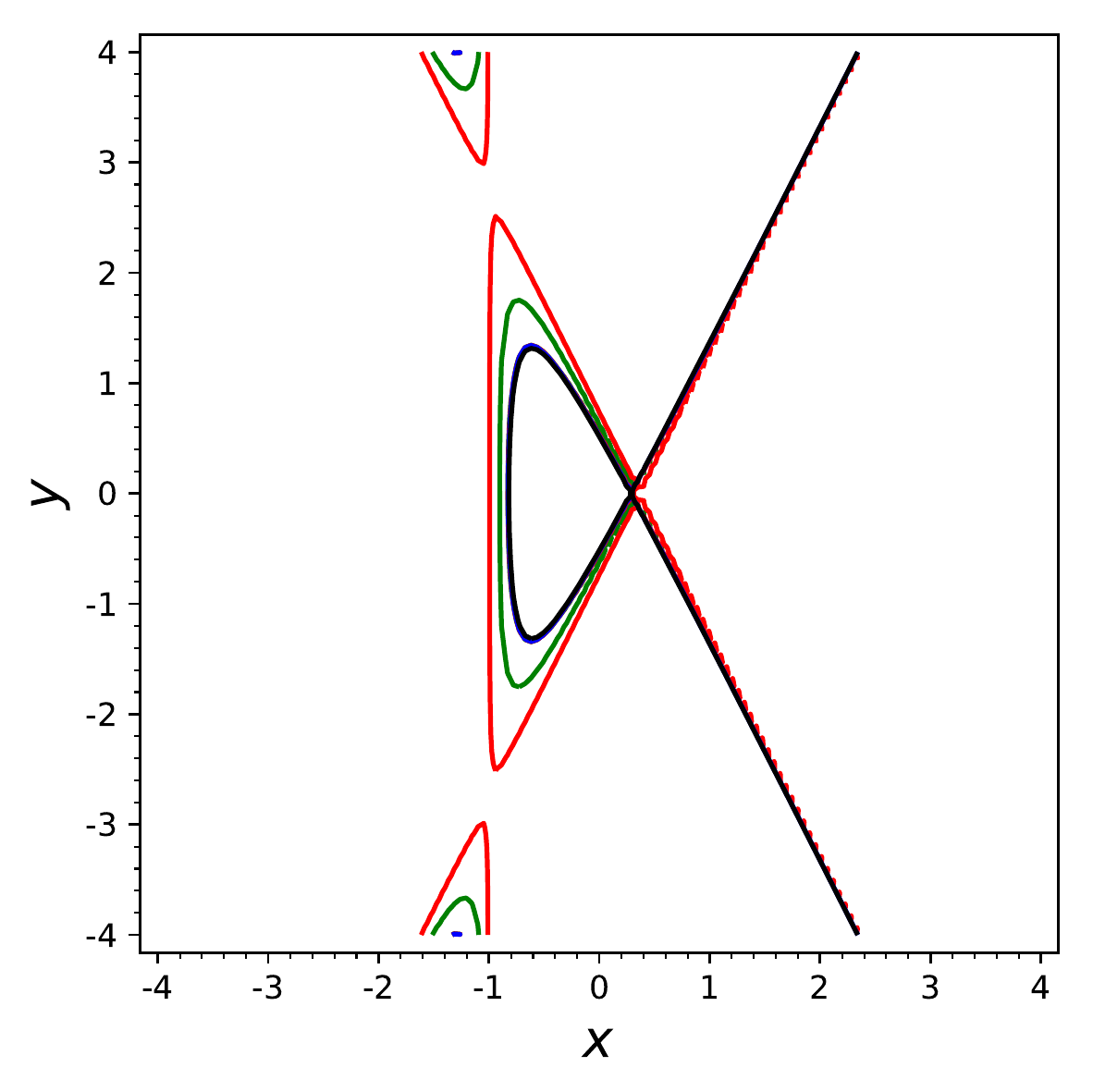}
\includegraphics[width=0.45\textwidth]{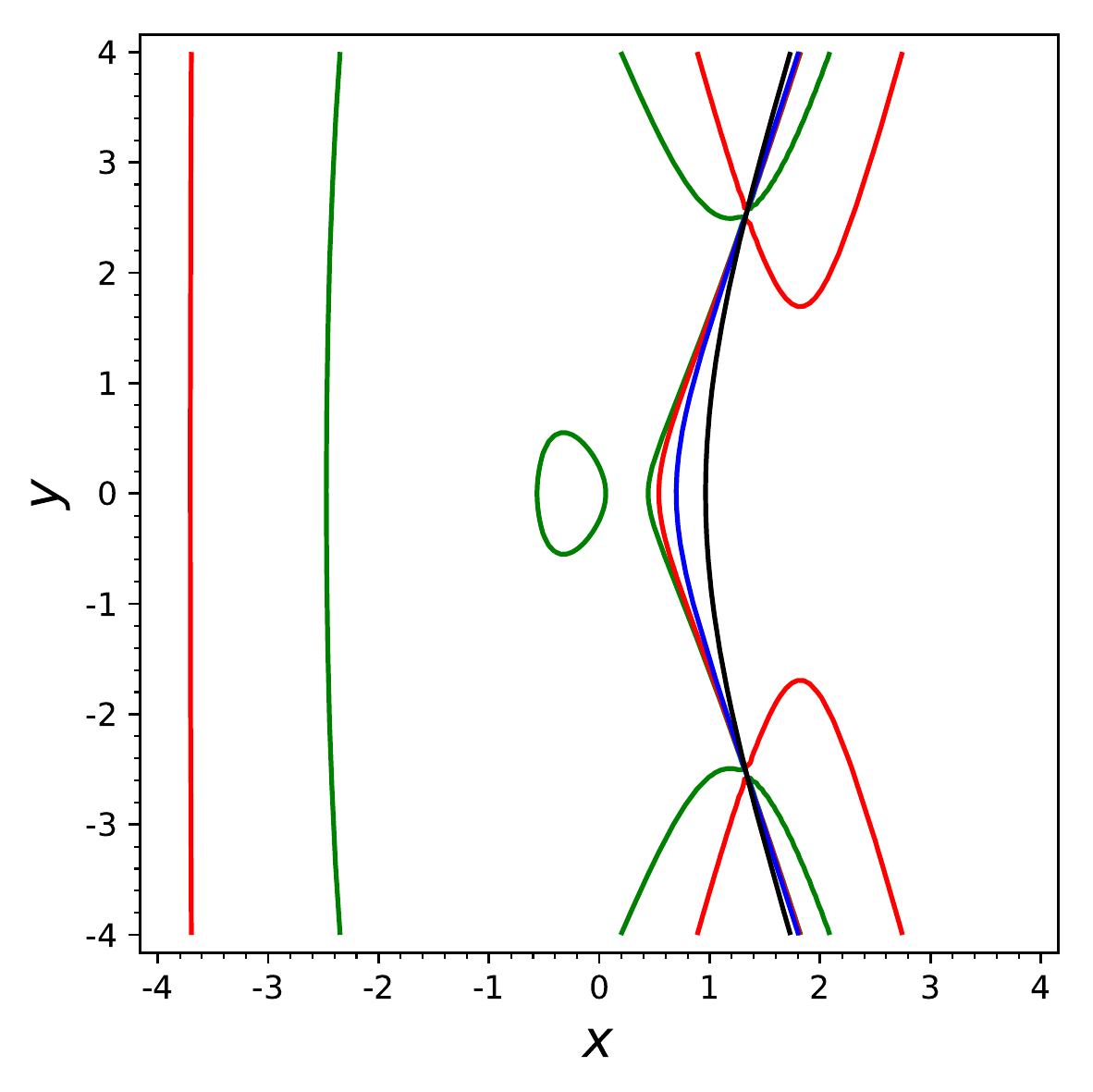}
\caption{Equiperiodic curves for Example \ref{ex:Fn} at the step $\dt=1$ (left) and at the step $\dt =0.5$ (right). 
%($F_5$ in black, $F_6$ in red, $F_7$ in blue, $F_8$ in green). 
}
\label{fig:wp:2}
\end{figure}

\begin{example}
\label{ex:Fn}
For the linear oscillator and $\wp$-oscillator, the equiperiodic set of the order  $2$ and $3$ is empty. At $n=3$ the curve equation degenerates into 
\[
3\dt^4 - 4=0,
\]
which agrees with the aforementioned circumstance: at $n=4$ the step is independent of the initial data (Example \ref{ex:wp:2}). At $n=5$ the equiperiodic set appears to be an elliptic curve
\begin{equation}
\label{eq:f5}
\begin{aligned}
& 27\dt^{10}x - 432\dt^8xy^2 + 432\dt^8x^2 + 1728\dt^6x^3 + 27\dt^8 \\
- &  432\dt^6y^2 - 936\dt^6x + 168dt^4 + 240\dt^2x - 80=0.
\end{aligned}
\end{equation}
This curve and the equiperiodic curves at $n=6,7$ and $8$ are plotted in Fig.~\ref{fig:wp:2}.  \alert{The degrees of $F_n \in \mathbb{Q}[\dt][x,y]$ are presented in Table \ref{t:Fn}. Due to the degree grows, the equiperiodic curves do not belong to the same sheaf, linear or irrational.}   
\end{example}

\begin{example}
In the case of the Jacobi oscillator, the equiperiodic set consists of some surface in the space $ pqr $, to which coordinate lines should be added. The expressions obtained are very cumbersome.
\end{example}

The very definition implies the simplest properties of equiperiodic sets.
\begin{itemize}
\item If some point of an approximate solution belongs to an equiperiodic set, then this solution has period $ n $ and all its points belong to this set. Therefore, equiperiodic sets are integral sets for an approximate solution. Thus, on periodic solutions, conservation laws are satisfied, but different from those known for the continuous model.
\item If equiperiodic sets of orders $ n'$ and $ n''$ intersect in a nonsingular point of the transformation, then the solution, which started at the intersection point, must have periods $ n' $ and $ n '' $ simultaneously. Therefore, the numbers $ n'$ and $ n''$ must have a common divisor $ n''' $, and the intersection point itself must lie on an equiperiodic set of order $ n '''$.
\end{itemize}

If we consider, e.g.,  the $ xy $ plane of a $ \wp $-oscillator with a fixed value of  step $ \dt $,  we will observe a countable number of equiperiodic curves $ F_n (x, y, \dt) = 0 $. The solutions, which started at a point on the $ F_n $ curve, run through exactly $ n $ points of this curve and come back. The solutions that do not fall on these curves will be aperiodic. 

\section{Discussion and conclusion}

In classical mechanics, attempts were periodically made to consider Newton's differential equations as difference equations, considering $ \dt $ as a small but finite quantity
\cite{Feynman}. 
However, if we replace the equation  $\tfrac{dx}{dt} = f(x)$
with $\hat x -x =f(x)\dt$,  all fundamental laws of nature are violated, including $ t $-symmetry and conservation laws. As a result, the difference model loses the well-known properties of the continuous model and one has to use a continuous model with the correct qualitative properties, and the discrete model is considered as an imperfect one suitable only for numerical calculations.  

The ultimate goal is to create discrete models that have the most important properties of mechanical models. These include undoubtedly the inheritance of algebraic conservation laws, $ t $-symmetry, reversibility and periodicity. As we found out earlier, it is impossible to combine invertibility and exact preservation of all algebraic integrals.
\cite{malykh-casc-2019}.
In this paper, as in the continuous case, we restricted the consideration to an integral manifold (Example \ref{ex:wp:1}) and considered difference schemes on the manifold. However, in the discrete case, one can go another way, rejecting the exact preservation of precisely those expressions that are preserved in the continuous case.

The starting point for this article was the observation that dynamical systems with a quadratic right-hand side can be appro\-ximated by reversible difference schemes with $ t $-symmetry. The approximate solutions found using these schemes are birational functions of the initial data over the entire phase space. This is surprising since
in the continuous case, for this property to appear, one had to restrict the phase space by an algebraic integral manifold (see Example \ref{ex:wp:1}).

For any initial data and any $ n \in \mathbb N $, computer algebra methods can find all possible values for the step $ \dt $ , at which the approximate solution is a sequence of points with a period $ n $. The experiments carried out demonstrate that for sufficiently large $ n $ the set of such steps is not empty and the minimum step tends to zero with increasing $ n $. Thus, for any initial data, one can specify a periodic approximate solution, arbitrarily close to the exact one in the uniform norm.

The sets of points in the phase plane, from which the solutions emerge, which are sequences with the period $ n $, were called equiperi\-odic in the text. Their study seems to us to be a very promising and beautiful task, which can clarify to what extent discretization according to a reversible difference scheme randomizes a completely integrable continuous problem. 

\begin{acknowledgments}
This work is supported by the Russian Science Foundation (grant no.~20-11-20257). 
\end{acknowledgments}

\bibliographystyle{ugost2008}
\bibliography{malykh}

\end{document}